\documentclass[12pt]{amsart}
\usepackage{amsbsy,amssymb,amscd,amsfonts,latexsym,amstext,delarray,
amsmath,epsfig}

\newtheorem{thm}{Theorem}[section]  
        
\newtheorem{cor}[thm]{Corollary}            
\newtheorem{lem}[thm]{Lemma}                    
\newtheorem{defn}[thm]{Definition}
\newtheorem{rem}[thm]{Remark} 
\newtheorem{ex}[thm]{Example}
\numberwithin{equation}{section} 

\def\R{{\mathbb R}}  

\def\Z{{\mathbb Z}} 
 
\def\C{{\mathbb C}}   
\def\N{{\mathbb N}}

\def\p{{\mathfrak p}}
\def\P{{\mathfrak P}}
\def\Ss{{\mathcal S}} 

\def\A{{\mathcal A}}
\def\G{{\mathcal G}}
\def\B{{\mathcal B}}
\def\Cl{{\mathfrak Cl}}

\DeclareMathOperator{\spinc}{{\rm Spin}^c}
\def\s{{\mathfrak s}}
\DeclareMathOperator{\Hom}{{\rm Hom}}

\DeclareMathOperator{\cl}{\bf c} 
\DeclareMathOperator{\Di}{{D\mkern -12mu /}_A} 
\DeclareMathOperator{\di}{{\partial}\mkern -10mu /} 

\DeclareMathOperator{\dist}{dist}
\DeclareMathOperator{\supp}{supp}

\begin{document}

\title[Weak UCP and perturbed monopole equations]{Weak UCP and
perturbed monopole equations}  
\author[B.~Booss--Bavnbek]{B.~Booss--Bavnbek}
\author[M.~Marcolli]{M.~Marcolli$^{\text{\dag}}$} 
\author[B.L.~Wang]{B.L.~Wang$^{\text{\ddag}}$}
\thanks{\noindent$^{\text{\dag}}$Partially supported by Humboldt
Foundation Sofja Kovalevskaja Award}
\thanks{$^{\text{\ddag}}$Supported by Australian Research
Council} 
\address{Bernhelm Booss--Bavnbek: Institut for matematik og fysik\\
Roskilde University, 4000 Ros\-kilde, Denmark}
\email{booss@mmf.ruc.dk}
\address{Matilde Marcolli: Max--Planck Institut f\"ur Mathematik \\
Vivatsgasse 7, 53111 Bonn, Germany}
\email{marcolli@mpim-bonn.mpg.de} 
\address{Bai--Ling Wang: Department of Pure Mathematics \\ 
University of Adelaide, Adelaide SA 5005, Australia}
\email{bwang@maths.adelaide.edu.au}
\maketitle

\begin{abstract} 
We give a simple proof of weak Unique Continuation
Property for perturbed Dirac operators, using the Carleman
inequality. We apply the result to a class of perturbations of the
Seiberg--Witten monopole equations that arise in Floer theory. 
\end{abstract} 

\section{Introduction}
We outline the content and motivation of the paper. The first part of
the paper presents a proof of weak Unique Continuation
Property for perturbed Dirac operators, based on classical methods
revolving around the Carleman inequality. In
the second part we give an application to the Morse--Smale--Witten
complex for the Chern--Simons--Dirac functional, whose homology
defines a gauge--theoretic invariant of 3--manifolds.

\subsection{UCP and Dirac operators} A linear or non-linear operator
$\mathfrak{D}$, acting on functions or sections of a
bundle over a compact or non-compact manifold $M$ has
the {\em weak Unique Continuation Property (UCP)} if
any solution $u$ of the equation $\mathfrak{D}u=0$ has
the following property: if $u$ vanishes on a non--empty open
subset $\Omega$ of $M$, it vanishes on the whole
connected component of $M$ containing $\Omega$.

There is also a notion of {\em strong UCP}, where, instead of assuming
that a solution $u$ vanishes on an open subset, one assumes
only that $u$ vanishes `of high order' at a point. 
The concepts of weak and strong UCP extend a fundamental property of
analytic functions to {\it some} elliptic equations other than the
Cauchy--Riemann equation, see also \S \ref{s-ucp} below.

Up to now, (almost) all work on UCP goes back to
two seminal papers \cite{Ca33}, \cite{Ca39} by
Torsten Carleman, establishing a Carleman--type
inequality (cf.~ our inequalities \eqref{e-8.1} and \eqref{e-8.1_pert} 
below). In this approach, the difference between weak and strong UCP
and the possible presence of more delicate non--linear perturbations
are related to different choices of the weight function in the
inequality, and to whether $L^2$--estimates suffice or $L^p$-- and 
$L^q$--estimates are required.

\medskip

There are different notions of operators of Dirac type. We shall not 
discuss the original hyperbolic Dirac operator (in the Minkowski
metric), but restrict ourselves to the elliptic case related to
Riemannian metrics. 

Recall that, if $(M,g)$ is a compact smooth Riemannian manifold
(with or without boundary) with $\dim M=m$, we denote by
$\Cl (M) =
\{ \Cl (TM_x,g_x) \}_{x\in M}$ the bundle of Clifford
algebras of the tangent spaces. For $E\to M$ a
smooth complex vector bundle of Clifford modules, the
{\it Clifford multiplication} is a bundle map $\cl:
\Cl (M) \to \Hom(E,E)$ which yields a representation
$\cl:
\Cl (TM_x,g_x) \to \Hom_{\C}(E_x,E_x)$ in each fiber. We
may assume that the bundle $E$ is equipped with a
Hermitian metric which makes the Clifford multiplication
skew--symmetric
\begin{equation}\label{e-skew-herm}
\langle \cl(v) s,s'\rangle = - \langle s,\cl(v)s'\rangle\qquad
\text{ for $v\in TM_x \text{ and } s\in E_x$}.
\end{equation}

Any choice of a smooth connection
\[
\nabla:{\rm C}^\infty (M;E) \to
{\rm C}^\infty(M;T^*M\otimes E)
\]
defines an {\em operator of Dirac type}  ${\mathfrak D}:=\cl\circ
\nabla $ under the Riemannian identification of the
bundles $TM$ and $T^*M$. 
In local coordinates we have
$
{\mathfrak D}:= \sum_{j=1}^m \cl(e_j) \nabla_{e_j}
$
for any orthonormal base $\{e_1,\dots,e_m\}$ of $TM_x$.
Actually, we may choose a local frame in such a way
that  
\[ \nabla_{e_j} = \frac{\partial}{\partial x_j} + \text{ zero
order terms} 
\]
for all $1\le j\le m$. So, locally, we have
\begin{equation}\label{e-dirac}
{\mathfrak D}:= \sum_{j=1}^m \cl(e_j) \frac{\partial}{\partial x_j} +
\text{ zero order terms }.
\end{equation}
It follows at once that the principal symbol
$\sigma_1({\mathfrak D})(x,\xi)$ is given by Clifford multiplication
with $i\xi$, so that any operator ${\mathfrak D}$ of Dirac type is
elliptic with symmetric principal symbol. Actually, if the connection
$\nabla$ is {\it compatible} with Clifford multiplication (i.e.
$\nabla\cl=0$), then the operator ${\mathfrak D}$ itself becomes
symmetric. We shall, however, admit non--compatible
metrics. Moreover, the {\it Dirac Laplacian} ${\mathfrak D}^2$ has
principal symbol $\sigma_2({\mathfrak D}^2)(x,\xi)$ given by the
Riemannian metric
$\| \xi \|^2$. So, it is scalar real (i.e.~a real
multiple of the identity) and elliptic.

\subsection{Motivation} Our main motivation for investigating weak UCP
for perturbed Dirac equations is an application to gauge theory of
3--manifolds. We outline briefly the context in which the question
arises.

Seiberg--Witten Floer homology is an invariant of 3--manifolds
defined as the homology of a Morse--Smale--Witten complex for the
Chern--Simons--Dirac functional, defined on an infinite dimensional
space of $U(1)$--connections and spinor sections. A detailed
construction of Seiberg--Witten Floer homology along with an analysis
of its main properties is given in \cite{MW}.  

Several interesting analytical problem are connected to the
construction of this invariant. One source of technical difficulties
is finding a suitable perturbation theory for the functional, in
order to have the Morse (or Morse--Bott) condition for the critical
points, and transversality of the spaces of flow lines. There are
different ways of treating this problem. 

One possible approach is the one followed in \cite{MW}. First one
realizes that it is fairly easy to achieve transversality for the set
of critical points of the Chern--Simons--Dirac functional (moduli
space of gauge classes of 3--dimensional Seiberg--Witten monopoles),
while it is more difficult to achieve transversality for moduli
spaces of flow lines. This observation leads to the idea of
perturbing the functional just enough to achieve transversality at
critical points, and then introducing a class of perturbation of the
4--dimensional Seiberg--Witten monopole equations, which does not come
from a perturbation of the functional. The perturbations constructed
this way have to satisfy a certain list of properties \cite[\S
2.3]{MW}, ensuring that the resulting moduli spaces of flow lines have
the desired transversality properties, so that the boundary operator
in the Morse--Smale--Witten complex can be defined by a counting
of flow lines.

In this approach, all the perturbations can be chosen so that, in the
Seiberg--Witten monopole equations, only the {\em curvature equation}
is perturbed, while leaving the {\em Dirac equation} unchanged.
A class of perturbations with these  properties
was introduced by Fr{\o}yshov in \cite{Fro}. 

Though this approach is effective in providing a working definition of
Floer homology and in the proof of topological invariance 
(in the equivariant case) \cite{MW}, it is very unnatural to use 
different perturbations for critical points and flow lines which do 
not come from a perturbation of the functional. Moreover, with the
perturbation defined by Fr{\o}yshov, it is hard to have good control of
the effect on the equations when the underlying 3--manifold is
modified, for instance by stretching a long cylinder, as in problems
related to cutting and pasting (surgery formulae).
In the interest of deriving formulae of this sort, it is better to
investigate other possible perturbation theories, which arise directly as
perturbations of the Chern--Simons--Dirac functional. A class of such
perturbation was proposed by Kronheimer in \cite{Kron}. Since these
affect both the curvature and the Dirac equation, and are both non--local
and non--linear, the question arises of how much delicate properties
of the Dirac operator, such as the Unique Continuation Property, may 
be affected by the presence of perturbation. In particular, since the 
weak UCP plays a role in all the transversality arguments, this seems 
an important technical point that needs to be understood.

This is the main point in our paper. We apply our elementary proof
of weak UCP for perturbed Dirac operators to a class of 
perturbations for Seiberg--Witten Floer theory that combines the
perturbations introduced in \cite{CW} and \cite{Kron}.  
  
Finally, it should be mentioned that a version of Seiberg--Witten
Floer theory that avoids the use of perturbations and deals directly
with the resulting excess intersection was developed recently by
Manolescu \cite{Manol}.  

\bigskip

\noindent{\bf Acknowledgment.} We thank Hubert Kalf for his generous
help and many extremely useful comments and suggestions. 

\section{Weak UCP for perturbed Dirac operators}

Clearly, not any arbitrary perturbation by a 0th order term  
preserves the {\em Unique Continuation Property (UCP)} 
as one can see by standard cases from first order ordinary 
differential equations, such as $\dot{u} -2\sqrt{|u|} =0$, or
$\dot{u} -3 u^{2/3}=0$. With a little more work, it is possible to
produce examples with non--uniqueness of {\em smooth} solutions.

In the affirmative, weak UCP can be proved in the following case.

\begin{thm}\label{t-main}
The weak UCP (i.e. UCP from open subsets) for 
solutions of a possibly tensored Dirac equation $\Di 
u=0$ on a smooth (not necessarily compact) manifold $M$ 
is preserved under perturbation with a non--linear 
and/or global term of 0th order of the form $\P(u)$ 
which can be estimated in the following way:
\begin{itemize}
\item  $\bigl|\P(u)(x) \bigr| \le 
P(u,x) |u(x)|$ with a real valued non--negative 
locally bounded function $P(u,\cdot)$ on $M$.
\end{itemize} 
\end{thm} 

\subsection{Guide to the Literature}
 
We review briefly three ways to prove Theorem \ref{t-main}. 
 
\subsubsection{} We can replace the Dirac operator by its square, 
the Dirac Laplacian. It has diagonal and real principal 
symbol. So we can apply the legendary parallel papers 
\cite{Ar57} and \cite{Co56} by Aronszajn and  
Cordes. Of course, in principle, a perturbation of a differential 
operator of first order by a 0th order term is 
essentially more delicate than a perturbation of a 
second order operator by 0th order. Fortunately, the 
two mentioned papers admit non--linear perturbations of 
first order of precisely that type which our 0th 
order perturbation of the Dirac operator will yield by 
squaring the perturbed operator.

Strictly speaking, Aronszajn and Cordes derived their 
results only for scalar equations and not for systems, 
but they mentioned that the proofs are similar for systems. 
Moreover, it should be mentioned that these methods are 
extremely hard because the goal of the two famous  
papers was the strong UCP (from a point). However, for 
the applications in Seiberg--Witten--Floer theory, weak UCP 
suffices. Finally, from a geometric point of view it 
seems inappropriate to destroy the beautiful geometric 
first--order structure of the Dirac operator by 
squaring. 

\subsubsection{} For weak UCP, an alternative and much  
simpler proof was given   
by Weck, \cite{We82} for any first order system of generalized 
Dirac type and a non--linear smooth perturbation satisfying the 
condition of Theorem \ref{t-main}. He transforms to a particular 
second order system and then establishes weak UCP for that system 
relatively easily.  
 
\subsubsection{} In the following we shall explain a third and 
completely direct proof method, where all 
arguments are carried out on the level of the perturbed 
Dirac operator, without resorting to second order 
operators.  
 
First we recall some features of the direct proof of 
weak UCP for solutions of the Dirac equation, given in 
\cite{BoWo93} and \cite{Bo00a}. Later on we explain the 
modifications necessary to cover also the non--linear 
case. 
 
\bigskip 

\subsection{The Unique Continuation Property - The 
Unperturbed Case}\label{s-ucp}

The weak UCP is one of the basic properties of an operator of Dirac 
type ${\Di}$. For $M=M_0\cup_{\Sigma}M_1$\,, it guarantees that 
there are no {\it ghost} solutions of ${\Di}u=0$, that is, 
there are no solutions which vanish on $M_0$ and have 
non--trivial support in the interior of $M_1$. This 
property is also called UCP {\it from open subsets} or 
{\it across any hypersurface}. 
For Euclidean (classical) Dirac operators 
(i.e., Dirac operators on $\R^m$ with constant coefficients
and without perturbation) the property follows by
squaring directly from the well--established UCP for the classical 
(constant coefficients and no potential) Laplacian.

In \cite[Chapter 8]{BoWo93} a very simple proof of 
the weak UCP for operators of Dirac type is given, inspired by 
\cite[Sections 6-7, in particular the proof of inequality  
(7.11)]{Ni73} and \cite[Section II.3]{Tr80}.  
We refer to \cite{Bo00a} for a further 
slight simplification and a broader perspective.  
The proof does not use advanced 
arguments of the Aronszajn/Cordes type regarding the 
diagonal and real form of the principal symbol of the 
Dirac Laplacian, but only the following product property 
of Dirac type operators (besides G{\aa}rding's inequality,  
see Remark \ref{r-gaarding}).

\begin{lem}\label{l-product} 
Let $\Sigma$ be a closed hypersurface of $M$ 
with orientable normal bundle. Let $t$ denote a normal 
variable with fixed orientation such that a bicollar 
neighborhood ${N}$ of $\Sigma$ is parameterized by $[- 
\epsilon,+\epsilon]\times \Sigma$. Then any operator of Dirac type can
be rewritten in the form 
\begin{equation}\label{e-product} 
{\Di}|_{{N}} 
=\cl(dt)\Bigl( \frac{\partial}{\partial t} + B_t + C_t \Bigr), 
\end{equation} 
where $B_t$ is a self--adjoint elliptic operator on the 
parallel hypersurface $\Sigma_t$, and $C_t: S|_{\Sigma_t} \to 
S|_{\Sigma_t}$ a skew--symmetric operator of $0$th order, 
actually a skew--symmetric bundle homomorphism. 
\end{lem}

\begin{rem}\label{r-homo-perturbation} {\em 
(a) It is worth 
mentioning that the product form \eqref{e-product} is invariant 
under perturbation by a bundle homomorphism. More precisely: 
Let $\mathfrak{D}$ be an operator on $M$ which can be written in 
the form \eqref{e-product} close to any closed hypersurface 
$\Sigma$\,, with $B_t$ and $C_t$ as explained in the preceding 
Lemma. Let $R$ be a bundle homomorphism. Then
\[
\bigl(\mathfrak{D}+R\bigr)|_N 
= \cl(dt)\Bigl(\frac{\partial}{\partial t} + B_t + 
C_t\Bigr) +\cl(dt)S|_N
\]
with $S|_N:=\cl(dt)^* R|_N$. Splitting $S=\frac 12(S+S^*) 
+\frac 12 (S-S^*)$ into a symmetric and a skew--symmetric part 
and adding these parts to $B_t$ and $C_t$, respectively, yields 
the wanted form of $\bigl(\mathfrak{D}+R\bigr)|_N$.

\noindent (b) For operators of Dirac type, it is well known that a 
perturbation by a bundle homomorphism is equivalent to modifying 
the underlying connection of the operator. This gives an 
alternative argument for the invariance of the form 
\eqref{e-product} for operators of Dirac type under perturbation by
a bundle homomorphism.

\noindent (c) By the preceding arguments (a), respectively (b), 
establishing weak UCP for sections belonging to the kernel of a 
Dirac type operator respectively an operator which can be written 
in the form \eqref{e-product} implies weak UCP for all 
eigensections. Warning: for {\em general linear elliptic} differential 
operators, weak UCP for ``zero-modes'' does not imply weak UCP for 
all eigensections.}
\end{rem}

To prove the weak UCP, in combination with the 
preceding lemma, the standard lines of the UCP 
literature can be radically simplified, namely with 
regard to the weight functions and the integration 
order of estimates. In Section \ref{s-ucp} we shall see that 
exactly these simplifications make it very easy to 
generalize the weak UCP to the perturbed case. 
We recall some decisive steps of the simplified proof 
in the unperturbed case. 

\smallskip

We want to show that, if $u\in{\rm C}^\infty(M,S)$ is a 
solution of ${\Di}u=0$, which vanishes on a non--empty open 
subset $\Omega$ of $M$, then it vanishes on the whole 
connected component of the manifold.

\smallskip

\subsubsection{} 
\begin{figure}
\epsfig{file=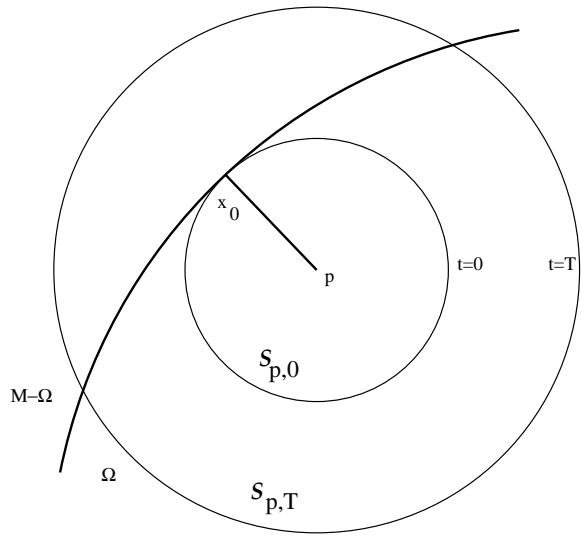}
\caption{\label{f-ucp} Local specification for the 
Carleman estimate}
\end{figure}

First we localize and convexify the situation and we 
introduce spherical coordinates (see Figure 
\ref{f-ucp}). Without loss of generality we may assume 
that $\Omega$ is maximal, namely the union of all open 
subsets where $u$ vanishes. If the solution $u$ does 
not vanish on the whole connected component containing 
$\Omega$, we consider a point 
$x_0\in \supp u\cap \partial\Omega$.  
We choose a point $p$ inside of $\Omega$ such that the ball 
around $p$ with radius $r:=\dist(x_0,p)$ is contained 
in $\overline{\Omega}$. We call the coordinate, running 
from $p$ to $x_0$\,, the {\it normal} coordinate and 
denote it by $t$. The boundary of the ball around $p$ 
of radius $r$ is a hypersphere and will be denoted by 
${\Ss}_{p,0}$. It goes through $x_0$ which has a normal 
coordinate $t=0$. 

Correspondingly, we have larger hyperspheres 
${\Ss}_{p,t}\subset M$ for $0\le t\le T$ with $T >0$ 
sufficiently {\it small}. In such a way we have parameterized 
an annular region 
${N}_T:=\{{\Ss}_{p,t}\}_{t\in[0,T]}$ around $p$ of 
width $T$ and inner radius $r$, ranging from the 
hypersphere 
${\Ss}_{p,0}$ which is contained in $\overline{\Omega}$, to 
the hypersphere ${\Ss}_{p,T}$ which cuts deeply into 
$\supp u$, if $\supp u$ is not empty.   

\subsubsection{} Next, we replace the solution $u$ by a cutoff 
\begin{equation}\label{e-bump} 
v(t,y):=\varphi(t)u(t,y)  
\end{equation}  
with a smooth bump function $\varphi$ with $\varphi(t)=1$ for 
$t\le 0.8 \,T$ and $\varphi(t)= 0$ for $t\ge 0.9\, T$. Then 
$\supp v$ is contained in ${N}_T$\,. More precisely, it 
is contained in the annular region ${N}_{0.9\, T}$\,. 
Moreover, $\supp ({\Di}v)$ is contained in the annular 
region $0.8\,T\le t \le 0.9 \,T$\,. 

We recall the following Lemma (cf.~Lemma 5 and Lemma 6 in 
\cite{Bo00a}). The weak UCP for the {\em unperturbed} Dirac 
operator $\Di$ will then follow immediately from this result.  

\begin{lem}\label{l-8.6} 
Let ${\Di}:{\rm C}^\infty (M,E)\to {\rm C}^\infty(M,E)$ be a linear
elliptic differential operator of order 1 which can be written 
on ${N}_T$ in the product form \eqref{e-product}. Let 
$v$ be a spinor section in the domain of ${\Di}$, such that $\supp(v)
\subset N_T$. 
 
\noindent {\em (a)}  
Then for $T$ sufficiently small there exists a constant 
$C$ such that the {\em Carleman} inequality (see 
\cite{Ca33}, \cite{Ca39} for the original form) 
\begin{equation}\label{e-8.1} 
R\int_{t=0}^T\int_{{\Ss}_{p,t}} e^{R(T-t)^2}\, 
|v(t,y)|^2\, dy\,dt 
\le C\int_{t=0}^T\int_{{\Ss}_{p,t}} e^{R(T-t)^2}\, 
|{\Di}v(t,y)|^2\, dy\,dt 
\end{equation} 
holds for any real $R$ sufficiently large. 
 
\noindent {\em (b)} Let $u$ be a solution of ${\Di} u =0$, with
$\supp (\varphi u) \subset N_{0.9T}$, for $\varphi$ the cutoff function
as in \eqref{e-bump}. If \eqref{e-8.1} holds for any sufficiently 
large $R>0$, then $u$ is equal $0$ on ${N}_{T/2}$\,. 
\end{lem} 

\bigskip 

Notice that, in the proof of \eqref{e-8.1} we do not assume $v$ is a 
solution, and we do not even assume it is smooth: the Carleman 
inequality is valid for {\em any} section with 
sufficiently small support, whether it is the cutoff of a
solution of the homogeneous equation or not. 
 
In the following section, we will return to some 
special features of the proof of the preceding lemma 
and show that the lemma remains true when an admissible 
(i.e. suitably bounded) perturbation is introduced.

\subsection{The Perturbed Case}\label{s-ucp_pert}  

We replace the equation $\Di u=0$ by  
\begin{equation}\label{e-dirac_pert} 
\widetilde \Di u := \Di u + \P_A (u) =0\,, 
\end{equation} 
where $\P_A$ is an {\em admissible} perturbation, in the following
sense.  

\begin{defn}\label{d-pert} 
A perturbation is {\em admissible} if it satisfies the 
following estimate: 
\begin{equation}\label{e-pert} 
\bigl|\P_A(u)|_x\bigr| \le P(u,x) |u(x)| \qquad 
\text{for $x\in M$} 
\end{equation} 
with a real--valued function $P(u,\cdot)$ which is 
locally bounded on $M$ for each fixed $u$.
\end{defn} 
 
\begin{ex}\label{ex-pert} {\em Some typical examples of
perturbations 
satisfying the admissibility condition of Definition \ref{d-pert} are:
\begin{itemize}
\item Consider a {\em non--linear perturbation} 
\[ 
\P_A(u)|_x := \omega(u(x))\cdot u(x)\,, 
\] 
where $\omega(u(x))|_{x\in M}$ is a (bounded) function 
which depends continuously on $u(x)$, for instance, for a 
fixed (bounded) spinor section $a(x)$ we can take 
\[ 
\omega(u(x)):=\langle u(x),a(x)\rangle 
\] 
with $\langle \cdot,\cdot\rangle$ denoting the Hermitian 
product in the fiber of the spinor bundle over the base 
point $x\in M$. This satisfies \eqref{e-pert}.
\item Another interesting example is provided by 
(linear) non--local perturbations with  
\[ 
\omega(u,x)=\left| \int k(x,z) u(z) dz\right| 
\] 
with suitable integration domain and integrability of 
the kernel $k$. These also satisfy \eqref{e-pert}.
\item Clearly, an unbounded perturbation may be both 
non--linear and global at the same time. This will, in fact, be the
case in our main application. In all these 
cases the only requirement is the estimate \eqref{e-pert} 
with bounded $\omega(u(\cdot))$. \end{itemize} } 
\end{ex} 
 
\medskip

We now show that (admissible) 
perturbed Dirac operators always satisfy the weak 
Unique Continuation Property. 

\begin{thm}\label{t-pert} 
Let $\Di$ be an operator of Dirac type and $\P_A$ an 
admissible perturbation. Then any solution $u$ of the 
perturbed equation \eqref{e-dirac_pert} vanishes 
identically on any connected component of the 
underlying manifold if it vanishes on a non--empty open subset of 
the connected component. 
\end{thm} 

\begin{proof} Let $u$ be a solution of the  perturbed 
equation which vanishes on an open non--empty set $\Omega$. 
We make the same construction as in the unperturbed 
case with a point $x_0\in \Omega\cap\supp(u)$, a point $p$ 
nearby in the interior of $\Omega$, a normal coordinate 
$t$, hyperspheres $\Ss_{p,t}$, a small positive number 
$T$, a cut--off function $\varphi$ with $\varphi(t)=1$ for 
$t\le 0.8 \,T$ and $\varphi(t)=0$ for $t\ge 
0.9 T$. Then we consider the cut--off solution 
$v=\varphi\cdot u$. Note that the support of $v$ is in the 
interior of the hypersphere $\Ss_{p,T}$ but $v$ is not 
a true solution. 

\medskip

We argue very much like in \cite[Proof of Lemma 6]{Bo00a}. To
begin with, we have

\begin{multline}\label{e-main1}
e^{R T^2/4}\, \int_0^{\frac T2}\int_{{\Ss}_{p,t}}
\|u(t,y)\|^2\, dy\,dt
= \int_0^{\frac T2}\int_{{\Ss}_{p,t}} e^{R T^2/4}\,
|u(t,y)|^2\, dy\,dt\\
\le \int_0^T\int_{S_{p,t}} e^{R(T-t)^2}\,
|\varphi u(t,y)|^2\, dy\,dt =: I
\end{multline}
We apply the Carleman type inequality of Lemma \ref{l-8.6} 
\begin{equation}\label{e-main2} 
I=\int_{t=0}^T\int_{{\Ss}_{p,t}} e^{R(T-t)^2}\, 
|\varphi u(t,y)|^2\, dy\,dt 
\le \frac CR\int_{t=0}^T\int_{{\Ss}_{p,t}} e^{R(T-t)^2}\, 
|{\Di}(\varphi u)(t,y)|^2\, dy\,dt .
\end{equation} 
We assume that $u$ is a solution of the perturbed equation
$\Di u +\mathfrak{P}_A(u) = 0$, hence
\[
\Di(\varphi u) = \varphi \Di u + \cl(dt)\varphi' u
= - \varphi \mathfrak{P}_A(u)  + \cl(dt)\varphi' u.
\]
Inserting in \eqref{e-main2} yields
\[
I \le \frac {2C}R\int_{t=0}^T\int_{{\Ss}_{p,t}} e^{R(T-t)^2}\, 
\bigl(|\varphi(t)\P_A(u)(t,y)|^2 +
|\cl(dt)\varphi'(t)u(t,y)|^2\bigr)\, dy\,dt .
\]
Now we exploit our assumption 
\begin{equation}\label{e-main3}
\bigl|\P_A(u)(x)\bigr| \leq P(u,x)|u(x)|\quad \text{for $x\in
M$}
\end{equation}
about the perturbation with locally bounded $P(u,\cdot)$, say
\[
|P(u, (t,y))|\leq C_0:=\max_{x\in K}|P(u,x)|
\quad\text{for all $y\in\Ss_{p,t},\, t\in[0,T]$}
\]
where $K$ is a suitable compact set. We obtain at once
\begin{align*}
\bigl(1-\frac{2CC_0}R\bigr)I 
&\leq \frac {2C}R\int_{t=0}^T\int_{{\Ss}_{p,t}} e^{R(T-t)^2}\, 
 |\cl(dt)\varphi'(t)u(t,y)|^2\, dy\,dt\\
&\leq \frac {2C}R e^{RT^2/25}\, \int_{t=0}^T\int_{{\Ss}_{p,t}} 
|\cl(dt)\varphi'(t)u(t,y)|^2\, dy\,dt.
\end{align*}
Here we use that $\varphi'(t)=0$ for $0\leq t\leq 0.8T$ so
that we can estimate the exponential and pull it in front of
the integral. Returning to \eqref{e-main1} yields
$$
\int_0^{\frac T2}\int_{S_t} \|u(t,y)\|^2\, dy\,dt \le 
\frac R{R-2CC_0}\frac {2C}{R} e^{-21R {T^2}/{100}} \int_0^T
\int_{\Ss_{p,t}} |\cl(dt)\varphi'(t)u(t,y)|^2\, dy\,dt,
$$
which gives the result as $R\to\infty$.
\end{proof}

\medskip

\begin{rem}\label{r-main}
{\em The preceding proof consists only of a short modification
of the usual way one obtains weak UCP for the unperturbed
equation from the Carleman estimate. The point of our proof
for the perturbed case is that the unmodified Carleman
estimate for the unperturbed operator suffices (Lemma
\ref{l-8.6}a). In the Appendix to this Note we give an
alternative proof of Theorem \ref{t-pert} by modifying Lemma
\ref{l-8.6} and establishing a new Carleman inequality.}
\end{rem}

\subsection{Application to the linearization of Seiberg-Witten
equations}

As an immediate application of Theorem \ref{t-pert}, we show that,
at any solution to the Seiberg-Witten equations, the linearization 
operator also enjoys the weak UCP. We only give an account of
such linearization in the 3-dimensional case, since the corresponding claim 
in the 4-dimensional case can be established by a similar argument.

Let $Y$ be a closed, oriented 3-manifold equipped  with a Riemannian metric
$g$ and a $\spinc$ structure. A $\spinc$ structure on $(Y, g)$
is a lift of the $SO(3)$--frame bundle on $Y$ to a $\spinc(3)$-bundle $P$.
Note that $\spinc(3) = U(2)$. The determinant homomorphism
$\spinc(3) \to U(1)$ determines a principal $U(1)$-bundle, whose 
corresponding complex line bundle is called the determinant line bundle
of the $\spinc$ structure $P$.

Given a $\spinc$ structure $P$, there is an associated $\spinc$ bundle
\[
 W= P \times _{U(2)} \C^2,
\]
which is a complex vector bundle and a module over the bundle
of Clifford algebras on $(Y, g)$. The  bundle
of Clifford algebras on $(Y, g)$ can be identified with
the bundle of exterior algebras on $T^*Y$, but with a different
algebra structure: the Clifford relation.

With  the Levi-Civita connection $\nabla$ on the cotangent bundle
$T^*Y$, a $U(1)$-connection $A$ on the determinant bundle $det (W)$
determines a connection
$\nabla_A$ on $W$ such that, for $v$ and $\psi$ sections of $T^*Y$ and
$W$ respectively,  
$\nabla_A$ satisfies
\[
\nabla_A (\cl(v)\psi) = \cl(\nabla v) \psi +\cl( v) \nabla_A (\psi).
\]
Then $\nabla_A$ is called a $\spinc$ connection on $W$.
Applying the Clifford multiplication, we can define a Dirac operator
$\di _A = \cl \circ \nabla_A: \Gamma (W) \to \Gamma (W) $.

The Seiberg-Witten equations on $(Y, g)$,
 for a pair $(A, \psi)$ consisting of a  $U(1)$ connection $A$ on the
determinant line bundle of $\s$ and a spinor section
$\psi$ of $W$, is given by 
\cite{KM} \cite{Fro} \cite{MW}:
\begin{equation}
\left\{ \begin{array}{l}
\di_A \psi =0, \\
*F_A = \sigma (\psi, \psi), 
\end{array}\right.
\label{SW-3d}
\end{equation}
where $\sigma (\psi, \psi)$ is an $\R$-bilinear form on
sections of spinor bundle, an imaginary valued 1-form
on $Y$ given by
\[
\sigma(\psi, \psi)= \displaystyle{\frac i2} 
\Im(\langle \cl(e_i)\psi, \psi \rangle )e^i \in \Omega^1(Y, i\R).\] Note that
these equations are gauge invariant under the action of
gauge group ${\mathcal G}=Map(Y, U(1))$:
$ u(A, \psi) = (A-2u^{-1}du, u\psi)$ for any $u\in {\mathcal G}=Map(Y, U(1))$.

Let $(A, \psi)$ be a smooth solution to the  Seiberg-Witten equations
(\ref{SW-3d}), then the spinor $\psi$ satisfies
\begin{equation}
|\psi|^2 \leq max_{y\in Y} \{ 0, -s (y)\}
\label{bound}
\end{equation}
where $s (y)$ is the scalar curvature for $(Y, g)$ \cite{KM}.

The linearization of the Seiberg-Witten equations (\ref{SW-3d})
at $(A, \psi)$, together with the linearization of the gauge action,
gives rise to the following linear elliptic system of equations
for $(\alpha, \phi)$ (a pair of an imaginary valued 1-form $\alpha$
and a spinor section):
\begin{equation}\label{3d-linear}
\left\{ \begin{array}{l}
d^*a + i\Im \langle \psi, \phi \rangle =0,\\
*d\alpha -\sigma (\psi, \phi),\\
\di _A\phi + \frac{1}{2} \cl(\alpha)\psi =0.
\end{array}\right.
\end{equation}

Then we have the following weak unique continuation
result for the linearization \eqref{3d-linear} of the Seiberg-Witten 
equations, whose proof follows from Theorem \ref{t-pert} and the 
pointwise bound in \ref{bound}.

\begin{cor}
If $(A,\psi)$ is a solution of the Seiberg-Witten equations (\ref{SW-3d}),
then as a solution to \eqref{3d-linear}, $(\alpha, \phi)$ 
satisfies the weak unique continuation property.
\end{cor}

\section{Perturbations of the Chern--Simons--Dirac functional}

We consider perturbations of the Chern--Simons--Dirac functional of
the form proposed in \cite{CW} and \cite{Kron}. The setup is as in the 
last paragraph: we have a closed compact connected oriented smooth 
3--manifold $Y$, with a fixed $\spinc$--structure $\s$, and with a 
choice of a Riemannian metric. We consider the configuration space 
$\A$ of pairs $(A,\psi)$ of 
a $U(1)$--connection and a spinor section, with the action of the gauge
group $\G$ as above. The spaces $\A$ and $Lie(\G)$ are completed
in suitable Sobolev norms (see e.g.~\cite{MW}). Moreover, on the
non--compact 4--manifold $Y\times \R$ we consider the
$\spinc$--structure obtained as pullback of $\s$ on $Y$, and the
cylindrical metric. The configuration space is given by the set
of {\em finite energy} pairs $({\mathbb A},\Psi)$ of a
$U(1)$--connection and a spinor on $Y\times \R$, acted upon by the
corresponding gauge group. The finite energy condition consists of the
property that, after a gauge transformation that kills the $dt$
component of ${\mathbb A}$, the resulting $(A(t),\psi(t))$ has 
time derivative in $L^2$. In this 
case the configuration space can be topologized by suitable weighted
Sobolev norms. Since we do not need the details here, we refer to
\cite{MW}.   

\subsection{Case I}
We consider the space of functions 
\begin{equation}
\label{UVspace}
\bigcup_{N\ge b_1, K>0}{\rm C}^{\infty}( \R^{N} , \R) \times {\rm
C}^{\infty}(\R^K, \R), 
\end{equation} 
where $b_1$ is the first Betti number of the 3-manifold $Y$.

In order to have the correct setup for transversality arguments, we
need a Banach space of parameters, hence we select as 
perturbation parameter space ${\mathcal P}$ a subspace of
\eqref{UVspace} of functions $(\p_1, \p_2)$ with finite Floer
$\epsilon$-norm. 

The Floer $\epsilon$-norm of $(\p_1, \p_2)$ is defined as follows:
choose $\underline{\epsilon} = (\epsilon_k)_{k\in \N}$ to
be a given sequence of positive real numbers, and set
\[
\| (\p_1, \p_2) \|_\epsilon =\sum_{k\ge 0} (\epsilon_k sup |\nabla ^k \p_1| +
\epsilon_k sup |\nabla ^k \p_2|).
\]
Following the argument of Lemma 5.1 of \cite{Flo}, the sequence 
$\underline{\epsilon}$ can be chosen   
such that ${\mathcal P}$ is a Banach space and 
${\mathcal P} \cap L^2_k$ is dense in $L^2_k$.
We require that $\p_1$ is invariant under the actions of
$H^1(Y, Z)$ on $\R^{b_1}$.

\subsubsection{}
For any 3-manifold $Y$, endowed with a Riemannian metric $g$ and a
fixed $\spinc$ structure $\s$, we choose a complete $L^2$--basis
$\{\nu_j\}_{j=1}^{\infty}$ of imaginary--valued 1--forms on $Y$.  We
also choose a complete $L^2$--basis $\{ \mu_j\}_{j=1}^{\infty}$ for the
co-closed (imaginary-valued) 1-forms on $Y$. Under the Hodge
decomposition, we have
\[
\Lambda^1_{L^2}(Y, i\R) = H^1(Y, i\R) \oplus Im (d^*) \oplus Im (d),
\]
and the $\{ \mu_j\}_{j=1}^{\infty}$ span the space $H^1(Y, i\R) \oplus
Im (d^*)$.  

\subsubsection{}
Fix a $U(1)$-connection $A_0$ on the determinant bundle $det (\s)$. To each 
co-closed 1-form $ \mu_j$ we associate a function on the
configuration space ${\mathcal A}$, defined as
\[
\tau_j (A, \psi) = \int_Y (A-A_0) \wedge * \mu_j.
\]
For simplicity, we assume that the $\{ \mu_j\}_{j=1}^{b_1}$ form a
basis of $ H^1(Y, i\R)$, and $[*\mu_j] =0$ for $j>b_1$. It is easy to
see that the following properties are satisfied:

(1) $\tau_j$ is invariant under gauge transformations for $j> b_1$;

(2) the map
\[
(\tau_1,\cdots,  \tau_{b_1}): {\mathcal A} \to \R^{b_1},
\]
is equivariant with respect to the action of 
gauge transformations $\lambda: Y \to U(1)$ on ${\mathcal A}$ and the
action of the corresponding $[\lambda] \in H^1(Y, \Z) \cong \Z^{b_1}$
on $\R^{b_1}$ as a translation by 
$$([\lambda^{-1}d\lambda] \cup [*\mu_1], [\lambda^{-1}d\lambda] \cup
[*\mu_2], \cdots, [\lambda^{-1}d\lambda] \cup [*\mu_{b_1}])\cap
[Y]. $$

\subsubsection{}
To each imaginary-valued 1-form $\nu_j$, we associate a function $\zeta_j$
on ${\mathcal A}$, defined as the quadratic form in the spinor $\psi$,
\[
\zeta_j (\psi, \psi)= 
\int_Y \langle \cl(\nu_j) \psi, \psi \rangle ,
\]
where $\langle\cdot, \cdot\rangle$ is the Hermitian metric on the
space of spinors. It is easy to 
see that $\zeta_j$ is gauge invariant and real-valued. 

\subsubsection{}
Now we choose any function $\p_1\in {\rm C}^{\infty} (\R^N, \R)$ and
$\p_2\in {\rm C}^{\infty} (\R^K, \R)$ (for 
$N\ge b_1$, $K>0$)  where $\p_1$ is invariant under the action of $H^1(Y, \Z)$
on $\R^{b_1} \subset \R^N$. We can define
a function on $\A/\G$ as 
\[
\p_1( \tau_1, \ldots, \tau_N) + \p_2( \zeta_1, \ldots, \zeta_K).
\]

Let $CSD$ be the Chern-Simons-Dirac functional 
\begin{equation}
CSD (A, \psi) = -\displaystyle{\frac {1}{2} \int_Y} (A-A_0) \wedge 
(F_A +F_{A_0} ) + \displaystyle{\int_Y} \langle \psi, \di_A
\psi\rangle,
\label{CSD}
\end{equation}
where $\di_A$ is the self adjoint Dirac operator on the compact
3--manifold $Y$ twisted with the $U(1)$--connection $A$.
We perturb this functional as
\begin{equation}
\widetilde{CSD} (A, \psi) = CSD (A, \psi) +
\p_1( \tau_1, \ldots, \tau_N)+ \p_2( \zeta_1, \ldots,
\zeta_K). \label{CSDpert} \end{equation} 

\subsubsection{}
For each $(\p_1,\p_2) \in {\mathcal P}$, the gradient of \eqref{CSDpert}, with
respect to the $L^2$--metric on ${\A}/{\G}$, is computed in the
following lemma \cite{CW}.    

\begin{lem} \label{lem-gradient}
Consider the perturbed functional \eqref{CSDpert} with
$$ (\p_1,\p_2) \in {\mathcal P} \cap \left( {\rm C}^{\infty} (\R^N ,
\R)\times {\rm C}^{\infty} (\R^K, \R) \, \right). $$
The $L^2$--gradient of \eqref{CSDpert} at $(A, \psi)$ is given by
\begin{equation}\label{gradient}
\nabla \widetilde{CSD}(A,\psi) = \left\{\begin{array}{l}
 *F_A -\sigma (\psi, \psi)- \sum_{j=1}^{N} \frac{\partial
\p_1}{\partial \tau_j}\mu_j \\[2mm]  
\di_A \psi- \sum_{j=1}^{K} \frac{\partial \p_2}{\partial
\zeta_j}\cl(\nu_j)\psi  \end{array}\right. 
\end{equation}
\end{lem}

\subsubsection{} There are two types of perturbed Dirac equations that
we obtain from Lemma \ref{lem-gradient}. In fact, the critical points
equations for the functional \eqref{CSDpert} provide 
perturbed 3--dimensional Seiberg--Witten equations of the form
\begin{equation} \label{3dP}
\left\{ \begin{array}{l}
*F_A = \sigma (\psi, \psi) +\displaystyle{
 \sum_{j=1}^{N} \frac{\partial \p_1}{\partial \tau_j}}\mu_j \\[2mm]
\di_A \psi= \displaystyle{\sum_{j=1}^{K} \frac{\partial \p_2}{\partial
\zeta_j}} \cl(\nu_j)\psi,
\end{array}\right.
\end{equation}
where $\sigma (\psi, \psi)$ is a quadratic form in the spinor given
by $\sum_i \langle \cl(e_i)\psi,\psi\rangle e^i$, 
in dual local basis $\{ e_i \}$ and $\{ e^i \}$ of $TY$ and $T^*Y$. 
Similarly, the gradient flow equations for the functional
\eqref{CSDpert} correspond to perturbed 4--dimensional Seiberg--Witten
equations (in a temporal gauge) of the form
\begin{equation} \label{4dP}
\left\{ \begin{array}{l} \frac{dA(t)}{dt} = -
*F_{A(t)} + \sigma (\psi(t), \psi(t)) +\displaystyle{
 \sum_{j=1}^{N} \frac{\partial \p_1}{\partial \tau_j}}\mu_j \\[2mm]
\frac{d\psi(t)}{dt} = - \di_{A(t)} \psi(t) +
\displaystyle{\sum_{j=1}^{K} \frac{\partial \p_2}{\partial \zeta_j}}
\cl(\nu_j)\psi(t). 
\end{array}\right.
\end{equation}

\subsubsection{} We can prove the following weak unique continuation
result for solutions of the perturbed Seiberg--Witten equations
\eqref{3dP} and \eqref{4dP}.

\begin{thm}\label{UCP3d4d}
If $(A,\psi)$ is a solution of \eqref{3dP}, and the spinor $\psi$
vanishes on a non-empty open set, then $\psi$ vanishes everywhere on the
compact connected 3--manifold $Y$. Suppose then that $({\mathbb
A},\Psi)$ is a pair of connection and spinor on the non--compact
four--dimensional manifold $Y\times \R$, with cylindrical metric
and spinor bundle $S=S^+\oplus S^-$ obtained by pulling back the
$\spinc$--structure $\s$ on $Y$. If $({\mathbb A},\Psi)$ is gauge
equivalent to $(A(t),\psi(t))$ in a 
temporal gauge, satisfying \eqref{4dP}, and $\Psi$ vanishes on a 
non-empty  open set, then $\Psi$ vanishes everywhere on $Y\times \R$.
\end{thm}

\begin{proof} The result follows in both cases directly from Theorem
\ref{t-pert}, with the non--linear and global perturbation
\begin{equation} \label{PA3d} \P_A(\psi) := -
\displaystyle{\sum_{j=1}^{K} \frac{\partial 
\p_2}{\partial \zeta_j}} \cl(\nu_j)\psi, \end{equation}
so that the Dirac equation in \eqref{3dP} becomes of the form
$\widetilde \di_A \psi =0$, with 
$$ \widetilde \di_A \psi = \di_A \psi + \P_A(\psi). $$
The Dirac equation in \eqref{4dP} can also be written
as $\widetilde \Di \Psi =0$, with the perturbation term
$\P_A(\Psi)$ which only differs from \eqref{PA3d} by a unitary operator
$\cl(dt)$ and  
$$ \widetilde \Di \Psi = \Di \Psi - \P_A(\Psi), $$
with $\Di$ the Dirac operator on $Y\times \R$ and $\P_A(\Psi)
= \displaystyle{\sum_{j=1}^{K} \frac{\partial
\p_2}{\partial \zeta_j}} \cl(\nu_j)\Psi$,
$$ \Di : \Gamma(Y\times \R, S^+) \to \Gamma(Y\times \R, S^-). $$

\end{proof}

\subsection{Case II} Here we consider a version of the perturbations
introduced in \cite{Kron}, where, in addition to the perturbations 
$\p_1(\tau_1,\ldots,\tau_N)$ and $\p_2(\zeta_1,\ldots,\zeta_K)$ 
defined as in \S 3.1.2 and 3.1.4, we introduce further
perturbations depending on both connection and spinor, described as
follows. 

\subsubsection{} \label{Green} Let $G$ be the Green operator for the
ordinary Laplacian on $Y$. If $\G_e$ is the identity component of the
gauge group $\G$, we can consider the subgroup ${\mathcal H} \subset
\G_e$  
$$ {\mathcal H}=\{ \lambda =e^{if}, f: Y \to \R \text{ such that }
\int_Y f =0 \}, $$ 
with $\G/{\mathcal H}\cong U(1) \times H^1(Y,\Z)$. Consider also a
fixed family of locally bounded spinors $\{\psi_i\}_{i=1, 2, \cdots }$ 
on $Y$. We choose the family $\{\psi_i\}_{i=1, 2, \cdots }$ so that
they form a complete $L^2$--basis of the space of spinor sections, and
they are eigenvectors of the fixed Dirac operator $\di_{A_0}$, where
$A_0$ is the fixed $U(1)$-connection that appears in the definition
of Chern-Simons-Dirac functional \eqref{CSD}. We set
\begin{equation}\label{eta}
 \eta_i (A,\psi):= \int \langle e^{-G d^* (A-A_0)} \psi_i , \,
\psi \rangle. \end{equation}

The following properties are satisfied:

(1) $\eta_i$ is invariant under ${\mathcal H}$;

(2) the map
$$ (\eta_1,\ldots, \eta_L): \A \to \C^L $$
is equivariant with respect to the action of $\G$ on $\A$ and the
corresponding action of $U(1)\times H^1(Y,\Z)$ on $\C^L$.

\subsubsection{} We choose a class of functions $\p_3: \C^L\to \R$ which
is invariant with respect to the action of $U(1)\times H^1(Y,\Z)$. 
We consider the corresponding perturbation term
$$ \p_3(\eta_1,\ldots, \eta_L) :\A \to  \R, $$
and the resulting perturbed Chern--Simons--Dirac
functional of the form 
\begin{equation}
\widetilde{CSD} (A, \psi) = CSD (A, \psi) +
\p_1( \tau_1, \ldots, \tau_N)+ \p_2( \zeta_1, \ldots,
\zeta_K) + \p_3 ( \eta_1,\ldots, \eta_L ). \label{CSDpert2} 
\end{equation}

\begin{lem} \label{gradient2lem}
The $L^2$--gradient of the perturbed functional \eqref{CSDpert2}, with
the additional perturbation $\p_3$, is given by
\begin{equation}\label{gradient2}
\nabla \widetilde{CSD}(A,\psi) = \left\{\begin{array}{l}
 *F_A -\sigma (\psi, \psi)- \sum_{j=1}^{N} \frac{\partial
\p_1}{\partial \tau_j}\mu_j - \sum_{\ell=1}^{L}\frac{\partial
\p_3}{\partial \eta_\ell }d G \left(i\Im \langle e^{-Gd^*
(A-A_0)}\psi_\ell, \psi\rangle \right)\\[2mm]
\di_A \psi- \sum_{i=1}^{K} \frac{\partial \p_2}{\partial
\zeta_i}\cl(\nu_i)\psi - \sum_{\ell=1}^{L}\frac{\partial
\p_3}{\partial \eta_\ell} e^{-Gd^* (A-A_0)} \psi_\ell
\end{array}\right. 
\end{equation}
\end{lem}

Thus, the Dirac equations derived from the perturbed critical point
equations 
$$\nabla\widetilde{CSD}(A,\psi) =0$$ 
and perturbed flow lines equations 
$$\frac{d}{dt}(A(t),\psi(t))+\nabla\widetilde{CSD}(A(t),\psi(t))
=0$$ 
for the functional \eqref{CSDpert2} are, respectively, of the form
\begin{equation}\label{3diracP}
\di_A \psi- \sum_{i=1}^{K} \frac{\partial \p_2}{\partial
\zeta_i}\cl(\nu_i)\psi - \sum_{\ell=1}^{L}\frac{\partial
\p_3}{\partial \eta_\ell} e^{-Gd^* (A-A_0)} \psi_\ell =0
\end{equation}
on the compact 3--manifold $Y$ and, on the non--compact four--manifold
$Y\times \R$,  
\begin{equation}\label{4diracP} 
\left(\frac{\partial}{\partial t} + \di_{A(t)} \right) \psi(t) 
-\sum_{i=1}^{K} \frac{\partial \p_2}{\partial 
\zeta_i}\cl(\nu_i)\psi(t) - \sum_{\ell=1}^{L}\frac{\partial
\p_3}{\partial \eta_\ell} e^{-Gd^* (A(t)-A_0)} \psi_\ell =0.
\end{equation}

\subsection{The UCP problem for Case II}

Unlike the perturbations $\p_1$ and $\p_2$ of Case I, the problem of
UCP is a lot more delicate for the perturbations $\p_3$ of Case II.

We begin with a very simple example, again taken from the theory of
ordinary differential equations, which shows that UCP may fail for a
toy model of the perturbations of Case II.
More precisely, we consider perturbations of the form
\begin{equation}\label{simplep3}
\mathfrak{P}_a(u)|_x := \langle u,a\rangle a(x), 
\end{equation}
for a fixed $L^2$--spinor $a(x)$, and with $\langle \cdot, \cdot
\rangle$ the $L^2$--inner product on spinors,
$$ \langle u,a\rangle =\int_Y \langle u(x), a(x) \rangle \, dv(x). $$
We take this class as a simplified version of \eqref{eta}, and of
the resulting $\p_3$. In any case, we are only looking at the Dirac
part of the Seiberg--Witten equations. 

\begin{ex} {\em
Let $a:[0,2]\to \R$ be a continuous function which
vanishes for $x\in [0,1]$ and satisfies
\[
\int_1^2 a(s)\,ds=\sqrt{2}\,.
\]
Then the function
\[
u(x):=\begin{cases}
0&\text{for $x\in [0,1)$}\\
\int_1^x a(s)\,ds&\text{for $x\in [1,2]$}
\end{cases}
\]
belongs to $C^1([0,2])$ and satisfies
\[
u'(x) = \langle u,a\rangle a(x)\qquad\text{(for $x\in
[0,2]$)},
\]
since
\[
\langle u,a\rangle = \int_0^2u(t)a(t)\,dt=\frac 12
u^2(2)=1\,.
\] }
\end{ex}

This shows that, already in the one--dimensional case and with the
simplified perturbations of the form \eqref{simplep3}, weak UCP
fails. This means that, in general, it will be difficult to expect
(even weak) UCP for the perturbations of Case II.

\medskip

However, there are special conditions under which perturbations of the
form \eqref{simplep3} do still satisfy weak UCP:

\begin{lem}\label{condition}
Assume that one of the following conditions is satisfied:
\begin{itemize}
\item the spinor $a$ does not vanish on any open subset, or 
\item The support of $a$ is contained in the interior of the 
support of the solution $u$ under consideration, or, alternatively 
put, there exists a positive constant $C_0$ such that 
$|a(x)| \leq C_0 |u(x)|$  for all $x\in M$.
\end{itemize}
Then the Dirac equation with a perturbation of the form
\eqref{simplep3} satisfies the weak UCP.
\end{lem}

\begin{proof} The second listed condition makes our perturbation 
admissible in the sense of Definition \ref{d-pert} (cf.~the second
item of Example \ref{ex-pert}), and weak UCP follows from Theorem
\ref{t-pert}. 

Now we deal with the first listed condition.
Let $u$ be a solution of the perturbed equation 
\[
\di_A u(x) + \langle u,a\rangle a(x)=0 ,\quad x\in M.
\]
If $\Omega$ is open in $M$ and $u$ vanishes identically on $\Omega$, 
then either $\langle u,a\rangle$ must vanish or $a$ must vanish 
identically on $\Omega$. In the first case we have no longer a 
perturbation and we are left with Lemma 2.3 which guarantees weak UCP.

The second case is excluded by the first listed condition of our 
lemma.
\end{proof}

By the choice of the $\{\psi_i\}_{i=1, 2,
\cdots}$, we know that each $\psi_i$, as an eigenvector of the
Dirac operator $\di_{A_0}$, satisfies the weak UCP by the argument of 
Remark \ref{r-homo-perturbation}, so  the additional perturbation 
in \eqref{3diracP} and \eqref{4diracP} satisfies the
first condition of Lemma \ref{condition}. The argument given in the
proof of Lemma \ref{condition} 
for perturbations of the form \eqref{simplep3}
extends easily to the more general case of perturbed equations of the
form \eqref{3diracP} and \eqref{4diracP}.

Thus, we have proved the following weak unique continuation
result for solutions of the perturbed Seiberg--Witten equations
\eqref{3diracP} and \eqref{4diracP}. 

\begin{thm}\label{UCP3d4dP}
If $(A,\psi)$ is a solution of \eqref{3diracP}, and the spinor $\psi$
vanishes on an open set, then $\psi$ vanishes everywhere on the
compact connected 3--manifold $Y$. Suppose then that $({\mathbb
A},\Psi)$ is pair of a connection and spinor on the non--compact
four--dimensional manifold $Y\times \R$, with cylindrical metric
and spinor bundle $S=S^+\oplus S^-$ obtained by pulling back the
$\spinc$--structure $\s$ on $Y$. If $({\mathbb A},\Psi)$ is gauge
equivalent to $(A(t),\psi(t))$ in a
temporal gauge, satisfying \eqref{4diracP}, and $\Psi$ vanishes on an open
set, then $\Psi$ vanishes everywhere on $Y\times \R$.
\end{thm}

\subsection{Linearizations}

A detailed discussion of the 
transversality results for monopole moduli spaces resulting from 
these perturbed equations will be presented elsewhere. However, in
view of such application, we state another result, which follows from
Theorem \ref{t-pert}.  
Recall that, in the proof of transversality for monopole equations
(see \cite[\S 2.3]{MW}), one argues that a possible element in the
Cokernel of the linearization must in fact vanish identically. By
varying the choice of the perturbation (transversality of the
universal moduli space), the question can be reduced to an argument
involving weak UCP for the adjoint of the (perturbed) Dirac operator
(we refer to  \cite[\S 2.3]{MW} for details). 

There is a conjecture by Laurent Schwartz \cite{Sc56} that weak 
UCP for an elliptic differential operator will always imply weak 
UCP for its formal adjoint. The matter was discussed in 
\cite{Bo65} in detail, but is still unsolved, see also 
\cite[Section 2.2]{Bo00a}. However, in the special case considered
here we can prove similar results for the linearizations and their
adjoints, applying Theorem \ref{t-pert}.

Consider the linearization of equations \eqref{3dP} and
\eqref{4dP}. This gives linear operators ${\mathfrak
d}_{(A,\psi)}$ and 
${\mathfrak D}_{({\mathbb A},\Psi)}$ respectively. Notice that the
linearization ${\mathfrak d}_{(A,\psi)}$ or its adjoint ${\mathfrak
d}_{(A,\psi)}^*$ acts on pairs of a 1--form and spinor
$(\alpha,\phi)$, and it involves a term of the form 
$\alpha \cdot \psi$. Thus, this operator is not of the type directly
considered in Theorem \ref{t-pert}, hence we shall not apply directly
the result to this operator. However, for the proof of transversality
(see e.g.~\cite{MW}) it is sufficient to know that the operators
${\mathfrak d}_{(A,\psi)}^*$ and ${\mathfrak D}_{({\mathbb A},\Psi)}^*$
acting on ``pure spinor elements'' $(0,\phi)$ satisfy the UCP
property, and this can be 
shown using Theorem \ref{t-pert}. The case of
${\mathfrak D}_{({\mathbb A},\Psi)}$ is completely analogous. 

\begin{cor}
Let the linear operators ${\mathfrak d}_{(A,\psi)}^*$ and 
${\mathfrak D}_{({\mathbb A},\Psi)}^*$ be the $L^2$--adjoint to
the linearizations of equation \eqref{3dP} and equation 
\eqref{4dP}, respectively. If an element $(0, \phi)$
of pure spinor part is in the kernel of ${\mathfrak d}_{(A,\psi)}^*$
or ${\mathfrak D}_{({\mathbb A},\Psi)}^*$, and the spinor $\phi$
vanishes on a non--empty open set, then $\phi$ vanishes identically on
$Y$ (or $Y\times \R$). 
\end{cor}

\begin{proof}
The adjoint operators $\di_A^*$ and $\Di^*$ are again of Dirac type. 
In particular the adjoint operators ${\mathfrak d}_{(A,\psi)}^*$ and 
${\mathfrak D}_{({\mathbb A},\Psi)}^*$, acting on elements of pure spinor 
part, are Dirac--type operators with admissible perturbations.
Then the weak UCP for an element in the kernel of ${\mathfrak
d}_{(A,\psi)}^*$ and ${\mathfrak D}_{({\mathbb A},\Psi)}^*$ follows from
Theorem \ref{t-pert}.

\end{proof}

\section{Appendix}

In this Appendix we give
an alternative proof of Theorem \ref{t-pert} which is also
instructive, though longer than the one we have given in the paper.

\smallskip
 
We prove that the perturbed Carleman 
inequality  
\begin{equation}\label{e-8.1_pert} 
R\int_{t=0}^T\int_{{\Ss}_{p,t}} e^{R(T-t)^2}\, 
|v(t,y)|^2\, dy\,dt 
\le C\int_{t=0}^T\int_{{\Ss}_{p,t}} e^{R(T-t)^2}\, 
|\widetilde{\Di}v(t,y)|^2\, dy\,dt 
\end{equation} 
holds for any real $R$ sufficiently {\it large}. 

In order to make the various steps of the proof more transparent, we
have subdivided the proof in short paragraphs. 

\subsubsection{}
First consider a few technical points. The Dirac 
operator $\Di$ has the form $G(t)(\partial_t+{\B}_t)$ on the 
annular region $\{\Ss_{p,t}\}_{t\in [0,T]}$, and it is 
obvious that we have wUCP for the operator 
$\partial_t+\B_t$ if and only if we have it for $\Di$. 
Moreover, we have by Lemma \ref{l-product} that 
$\B_t=B_t+C_t$ with a self--adjoint elliptic 
differential operator $B_t$ and an anti--symmetric 
operator $C_t$ of order zero, both on 
$\Ss_{p,t}$. 

\smallskip
 
\subsubsection{}Now make the substitution 
$ 
v=:e^{-R(T-t)^2/2}v_0 
$ 
which replaces \eqref{e-8.1_pert} by 
\begin{multline}\label{e-8.5} 
R\int_0^T\int_{\Ss_{t}} |v_0(t,y)|^2\, dy\,dt\\ \le 
C\int_0^T\int_{\Ss_{t}} \left|\frac {\partial v_0}{\partial t} 
+ \B_t v_0 + R(T-t)v_0+ e^{\frac R2(T-t)^2} 
\P_A(v)\right|^2\, dy\,dt\,. 
\end{multline} 
We shall denote the integral on the left hand side by $J_0$ 
and the integral on the right hand side by $J_1$. 

\subsubsection{}
Now we prove \eqref{e-8.5}. Decompose $\frac 
{\partial}{\partial t}+\B_t+R(T-t)+\P_A$ into its symmetric 
part $B_t+R(T-t)+\P_A$ and  
anti-symmetric part $\partial_t+C_t$. This gives 
\begin{align*} 
J_1&=\int\int  
\left|\left( 
\frac {\partial v_0}{\partial t} + C_tv_0 
          \right) +  
\left(B_tv_0 + R(T-t)v_0 
+ e^{\frac R2(T-t)^2} \P_A(v) \right) 
\right|^2\, dy\,dt\\ 
&= \underbrace{\int\int  
\left| 
\frac {\partial v_0}{\partial t} + C_tv_0 
\right|^2\, dy\,dt}_{J_{\operatorname{skew}}}  
+ \underbrace{\int\int  
| \left(B_t + R(T-t) \right)v_0}_{J_{\operatorname{sym}}}\\ 
&\qquad + e^{\frac R2(T-t)^2} \P_A(v)  |^2\, dy\,dt\\ 
&\qquad+ 2\Re\int\int  
\left\langle 
\frac {\partial v_0}{\partial t} + C_tv_0, \,\, B_tv_0 + R(T-t)v_0 
\right\rangle 
+ e^{\frac R2(T-t)^2} \P_A(v) \, dy\,dt\\ 
&=: J_{\operatorname{skew}} + J_{\operatorname{sym}} + 
J_{\operatorname{mix}}\,. 
\end{align*} 

\subsubsection{}
First we investigate the term $J_{\operatorname{mix}}$. 
We integrate by parts and use the identities for the 
real part $$ 
\Re\langle f,Pf\rangle =\frac 12\langle f,(P+P^*)f\rangle\,, 
$$  
and 
\[ 
2\Re\langle a,b\rangle \ge -2|\langle a,b\rangle| \ge -\left(\epsilon 
|b|^2 + \frac {|a|^2}{\epsilon}\right). 
\] 
We get, almost like in \cite[p.30]{Bo00a},  
\begin{equation}\label{e-jot_mix} 
J_{\operatorname{mix}}  = RJ_0 + J_3 + 
J_{\operatorname{skew.pert}}. 
\end{equation} 

\subsubsection{}
The first term on the right in \eqref{e-jot_mix} is 
\[ 
RJ_0 = R\int_0^T|v_0|^2\,dt. 
\] 
The second term on the right in \eqref{e-jot_mix} 
is the most delicate and 
was thoroughly analyzed in the unperturbed case: 
\[ 
J_3 := \int\int \left\langle v_0, \,\, -\frac {\partial B_t}{\partial 
t}v_0 + [B_t,C_t]v_0 \right\rangle \, dy\,dt\,. 
\] 
The third term on the right in \eqref{e-jot_mix} is 
\begin{align*} 
J_{\operatorname{skew.pert}} &-= 
2\Re\int\int \langle \frac {\partial v_0}{\partial t} + 
 C_tv_0 , \,\, e^{\frac R2(T-t)^2} \P_A(v) \rangle \,dy\, dt\\ 
&\ge -\int\int \left| \frac {\partial v_0}{\partial t} + 
 C_tv_0\right|^2 \,dy\, dt - \int\int\left| e^{\frac 
R2(T-t)^2} \P_A(v) \right|^2 \,dy\, dt. 
\end{align*}

\subsubsection{}
As a matter of fact the term $J_{\operatorname{skew}}$ 
was not used at all for establishing the Carleman 
inequality in the unperturbed case. Hence we can 
balance the first negative term in the preceding 
estimate for $J_{\operatorname{skew.pert}}$ with 
$J_{\operatorname{skew}}$. 

\smallskip

To balance the second negative term we decompose 
\begin{multline*} 
J_{\operatorname{sym}} 
= \epsilon \int\int \left|\left(B_t + 
R(T-t)\right)v_0\right|^2 \, dy\,dt\\ 
+ (1-\epsilon)\int\int \left|\left(B_t + 
R(T-t)\right)v_0\right|^2 \, dy\,dt\\ 
+ \int\int \left|e^{\frac R2(T-t)^2} \P_A(v)  
\right|^2\, dy\,dt \ +\ J_{\operatorname{sym.pert}} 
\end{multline*} 
with 
\begin{align*} 
J_{\operatorname{sym.pert}} 
&= 2\Re\int\int \langle \left(B_t + R(T-t)\right)v_0,\, 
 e^{\frac R2(T-t)^2} \P_A(v)  \rangle \, dy\,dt\\ 
&\ge -\epsilon\int\int \left|\left(B_t + 
R(T-t)\right)v_0\right|^2 \, dy\,dt\\ 
&\qquad - \frac 1{\epsilon} \int \int \left| e^{\frac 
R2(T-t)^2} \P_A(v) \right|^2 \,dy\, dt\,. 
\end{align*} 

\subsubsection{} Now the second (negative) term in the estimate for  
$J_{\operatorname{skew.pert}}$ balances exactly with 
the third term in the decomposition of 
$J_{\operatorname{sym}}$. Moreover, the first 
(negative) term in the estimate for 
$J_{\operatorname{sym.pert}}$ balances exactly with the 
first term in the decomposition of 
$J_{\operatorname{sym}}$. Note, that a fraction of 
$\int\int \left|\left(B_t + R(T- 
t)\right)v_0\right|^2 \, dy\,dt$ (say, the second term 
in the decomposition of $J_{\operatorname{sym}}$) was 
needed for balancing $J_3$ in the unperturbed case 
(see also Remark \ref{r-gaarding}). So, 
we must choose $\epsilon$ sufficiently small, say $\epsilon=\frac  
14$\,.

To balance the second (negative) term in the estimate 
for $J_{\operatorname{sym.pert}}$ we are left with the 
only positive term not yet (fully) used: $RJ_0$. 

\subsubsection{}
Let us first consider the case of admissible perturbations satisfying 
\eqref{e-pert}. Introducing the perturbation induces a final error   
\[ 
J_{\operatorname{err}} = 
\int_0^T\int_{\Ss_{p,t}} |v_0|^2 \left(R-\frac 
1{\epsilon}|P(v)|^2\right) \, dy\,dt 
\] 
with fixed $\epsilon$. For fixed $v$ we can choose a positive 
$R_0$ such that the error term $J_{\operatorname{err}}$ 
becomes $\ge 0$ for all $R\ge R_0$.

\medskip

\begin{rem}\label{r-gaarding} {\em 
In the {\em unperturbed} case, the completely elementary character  
of the arguments was made possible by exploiting one single  
non--elementary ingredient, G{\r a}rding's inequality (here the 
equivalence of the first Sobolev norm to the graph norm  
of any linear first order elliptic differential operator). It was used  
for estimating the term $J_3$ by $\frac 12 J_{\operatorname{sym}}$ 
and a fraction of $J_0$\,. For the {\em perturbed} case,  
it is worth mentioning that we do not need any kind of non--linear 
elliptic estimate but can keep the necessary modifications on the  
elementary level. }
\end{rem}

\end{document}